\documentclass{article}
\pdfoutput=1
\usepackage{geometry}     
\usepackage[round]{natbib}

\usepackage[utf8]{inputenc} %
\usepackage[T1]{fontenc}    %
\usepackage{hyperref}       %
\usepackage{url}            %
\usepackage{bbm, bm,color}
\usepackage{mathdots}
\usepackage{booktabs}       %
\usepackage{amsfonts}       %
\usepackage{nicefrac}       %
\usepackage{amsmath,amsthm}
\usepackage{graphicx}
\usepackage{caption}
\usepackage{subcaption}
\usepackage{enumitem}
\usepackage{algcompatible} %
\usepackage{wrapfig}
\usepackage{multicol}
\usepackage{algorithm}

\usepackage{color}
\hypersetup{ %
    pdftitle={},
    pdfauthor={},
    pdfkeywords={},
    pdfborder=0 0 0,
	pdfpagemode=UseNone,
    colorlinks=true,
    linkcolor=blue, %
    citecolor=blue, %
    filecolor=blue, %
    urlcolor=blue, %
    pdfview=FitH,
    pdfauthor={Simon Lacoste-Julien},
} 

\newlength{\Sfloatsep}
\newlength{\Stextfloatsep}
\newlength{\Sintextsep}
\newcommand{\R}{\mathbb{R}}

\newcommand{\X}{\mathcal{M}}

\newcommand{\prodscal}[2]{\left\langle#1,#2\right\rangle}
\newcommand{\x}{\bm{x}}
\newcommand{\y}{\bm{y}}

\newcommand{\s}{\bm{s}}
\newcommand{\xt}{\bm{x}^{(t)}}
\newcommand{\xtt}{\bm {x}^{(t+1)}}

\newcommand{\st}{\bm{s}^{(t)}}

\newcommand{\gap}{g_{t}}

\newtheorem*{rep@theorem}{\rep@title}
\newcommand{\newreptheorem}[2]{%
\newenvironment{rep#1}[1]{%
 \def\rep@title{#2 \ref{##1}}%
 \begin{rep@theorem}}%
 {\end{rep@theorem}}}

\newreptheorem{lemma}{Lemma'}
\newreptheorem{proposition}{Proposition'}
\newreptheorem{theorem}{Theorem'}
\newreptheorem{remark}{Remark'}

\newtheorem{definition}{Definition}
\newtheorem{theorem}[definition]{Theorem}

\DeclareMathOperator{\diam}{diam}

\DeclareMathOperator*{\argmin}{\arg\min}

\providecommand{\norm}[1]{\left\lVert#1\right\rVert}

\newcommand{\domain}{\mathcal{M}} %
\newcommand{\stepsize}{\gamma}

\newcommand{\Cf}{C_{\hspace{-0.08em}f}}

\newcommand{\dd}{\bm{d}}

\renewcommand{\r}{\bm{r}}

\newcommand{\innerProd}[2]{\left\langle #1 , #2 \right\rangle}

\newcommand{\indic}{\mathbbm{1}} %

\newcommand{\ignore}[1]{}%

\newcommand{\remove}[1]{} %

\title{Convergence Rate of Frank-Wolfe for Non-Convex Objectives}

\date{June 28, 2016}

\author{
Simon Lacoste-Julien \\
INRIA - SIERRA team\\
ENS, Paris 
}

\begin{document}
\maketitle

\begin{abstract}
We give a simple proof that the Frank-Wolfe algorithm obtains a stationary point at a rate of $O(1/\sqrt{t})$ on non-convex objectives with a Lipschitz continuous gradient. Our analysis is affine invariant and is the first, to the best of our knowledge, giving a similar rate to what was already proven for projected gradient methods (though on slightly different measures of stationarity).
\end{abstract}

\section{Introduction}
We consider the optimization problem:
\begin{equation} \label{eq:problem}
\min_{\x \in \domain} f(\x)
\end{equation}
where $f:\R^d\rightarrow\R$ is a continuously differentiable function over the domain $\domain$ that is \emph{convex} and \emph{compact}, but $f$ is potentially \emph{non-convex}. The Frank-Wolfe (FW) optimization algorithm proposed by~\citet{Frank:1956vp} (also known as conditional gradient method~\citep{demyanov1970approximate}), is a popular first-order method to solve~\eqref{eq:problem} while only requiring access to a \emph{linear minimization oracle} over~$\domain$, i.e., the ability to compute efficiently $\text{LMO}(\r) := \argmin_{\s \in \domain} \innerProd{\s}{\r}$.\footnote{Ties can be broken arbitrarily in this paper.} It has recently enjoyed a surge in popularity thanks to its ability to cheaply exploit the structured constraint sets $\domain$ appearing in machine learning applications, see~\citet{jaggi2013revisiting,lacoste2015global} and references therein. See also~\citet{lan2013FW} for a related survey.

We give the Frank-Wolfe algorithm with adaptive step sizes in Algorithm~\ref{alg:FW} (either with line-search or with a step size that minimizes an affine invariant quadratic upper bound). As~$\domain$ is convex, the iterates~$\xt$ stay in the feasible set~$\domain$ during the algorithm. For a convex function~$f$ with Lipschitz continuous gradient, the FW algorithm obtains a global suboptimality smaller than $\frac{2 C}{t+2}$ after~$t$ iterations~\citep[Theorem~1]{jaggi2013revisiting}, where $C \geq \Cf$ is the constant used in Algorithm~\ref{alg:FW} for the adaptive step size, and $\Cf$ is called the \emph{curvature constant} of~$f$ (defined in~\eqref{eq:Cf} below). On the other hand, we are not aware of any rates proven for Algorithm~\ref{alg:FW} in the case where $f$ is \emph{non-convex}. Examples of recent applications where the FW algorithm is run on a non-convex objective include multiple sequence alignment~\citep[Appendix B]{Alayrac16unsupervised} and multi-object tracking~\citep[Section~5.1]{chari15pairwise}. To talk about rates in the non-convex setting, we need to define a measure of non-stationarity for our iterates.

Consider the ``Frank-Wolfe gap'' of~$f$ at~$\xt$:
\begin{equation} \label{eq:FWgap}
g_t := \max_{\s \in \domain} \prodscal{\s-\xt}{-\nabla f(\x^{(t)})} .
\end{equation}
This quantity is a standard one appearing in the analysis of FW algorithms, and is computed for free during the FW algorithm (see Line~\ref{line:FWgap} in Algorithm~\ref{alg:FW}).
A point $\xt$ is a stationary point for the constrained optimization problem~\eqref{eq:problem} if and only if $g_t = 0$. Moreover, we always have $g_t \geq 0$. The FW gap is thus a meaningful measure of non-stationarity, generalizing the more standard $\| \nabla f(\xt) \|$ that is used for unconstrained optimization. An appealing property of the FW gap is that it is \emph{affine invariant}~\citep{jaggi2013revisiting}, that is, it is invariant to an affine transformation of the domain~$\domain$ in problem~\eqref{eq:problem} and is not tied to any specific choice of norm, unlike the criterion~$\| \nabla f(\xt) \|$. As the FW algorithm is also affine invariant~\citep{jaggi2013revisiting}, it is important that we state our convergence results in term of affine invariant quantities. In this paper, we show in Theorem~\ref{thm:gConvergence} below that the minimal FW gap encountered during the FW algorithm is $O(1/\sqrt{t})$ after $t$ iterations, that is:
\begin{equation}
\min_{0 \leq k \leq t} g_k  \leq \frac{\max\{2 h_0, C\} }{\sqrt{t+1}}  \quad \text{for $t \geq 0$} \, ,
\end{equation}
where $h_0 := f(\x^{(0)}) - \min_{\x \in \domain} f(\x)$ is the initial global suboptimality.

Another nice property of the FW gap $g_t$ is the following \emph{local suboptimality} property. If $\xt$ lies in a convex subset $\domain' \subseteq \domain$ on which $f$ is convex, then $g_t$ upper bounds the suboptimality with respect to the constrained minimum on~$\domain'$, that is, $g_t \geq f(\xt) - \min_{\x \in \domain'} f(\x)$ (by convexity).

\begin{figure}
\centering
\begin{minipage}[t]{.7\textwidth}
  \null
  \begin{algorithm}[H]
    \caption{Frank-Wolfe algorithm (with adaptive step sizes)}\label{alg:FW}
    \begin{algorithmic}[1]
      \STATE Let $\x^{(0)} \in \X$
      \FOR{$t=0 \ldots T$} 
      \STATE Compute $\st := \displaystyle\argmin_{\s \in \X} \prodscal{\s}{\nabla f(\x^{(t)})}$
      \STATE Let $\dd_t := \st - \xt$ \COMMENT{FW update direction}
      \STATE Compute $g_t := \left\langle \dd_t, - \nabla f(\x^{(t)}) \right\rangle$ \COMMENT{FW gap} \label{line:FWgap}
      \STATE \textbf{if} $ \gap \leq \epsilon$ \textbf{then} \textbf{return} $\x^{(t)}$ \COMMENT{approximate stationary pt.}
      \STATE Option I: Line-search: $\stepsize_t \in \displaystyle\argmin_{\stepsize \in [0,1]} \textstyle  f\left(\x^{(t)} + \stepsize \dd_t\right)$
      \STATE Option II: Set $\stepsize_t := \min\{ \frac{g_t}{C}, 1 \} \quad$ for some $C \geq \Cf$ given in~\eqref{eq:Cf}
      \STATE Update $\xtt := \xt + \stepsize_t \dd_t$
      \ENDFOR
      \STATE \textbf{return} $\x^{(T)}$
    \end{algorithmic}
  \end{algorithm}
\end{minipage}%
\end{figure}

\vspace{-3mm}
\section{Result}
\vspace{-2mm}
Before stating our convergence result, we review the usual affine invariant constant appearing in the convergence rates for FW methods. The \emph{curvature constant} $C_{f}$ of a continuously differentiable
function $f:\R^d\rightarrow\R$, with respect to a compact domain $\domain$, is defined as:
\begin{equation}\label{eq:Cf}
  \Cf := \sup_{\substack{\x,\s\in \domain,  ~\stepsize\in[0,1],\\
                      \y = \x+\stepsize(\s-\x)}} \textstyle
           \frac{2}{\stepsize^2}\big( f(\y)-f(\x)-\innerProd{\nabla f(\x)}{\y-\x}\big) \ .
\end{equation}
The assumption of bounded curvature $\Cf$ closely corresponds to a Lipschitz
assumption on the gradient of~$f$. %
More precisely, if~$\nabla f$ is $L$-Lipschitz continuous on $\domain$ with
respect to some arbitrary chosen norm $\norm{.}$ in dual pairing, i.e.
$\norm{\nabla f(\x) - \nabla f(\y)}_* \leq L \norm{\x-\y}$, then
\begin{equation} \label{eq:CfBound}
\Cf \le L \diam_{\norm{.}}(\domain)^2  \ ,
\end{equation}
where $\diam_{\norm{.}}(.)$ denotes the $\norm{.}$-diameter, see \citep[Lemma
7]{jaggi2013revisiting}. These quantities were normally defined in the context of convex optimization, but these bounds did not use convexity anywhere.

\begin{theorem}[Convergence of FW on non-convex objectives] \label{thm:gConvergence}
Consider the problem~\eqref{eq:problem} where $f$ is a continuously differentiable function that is potentially non-convex, but has a finite curvature constant~$\Cf$ as defined by~\eqref{eq:Cf} over the \emph{compact} convex domain~$\domain$. Consider running the Frank-Wolfe algorithm~\ref{alg:FW} with line-search (option I; then take $C := \Cf$ below) or with the step size that minimizes a quadratic upper bound (option II), for any $C \geq \Cf$. Then the minimal FW gap $\tilde{g}_t := \displaystyle \min_{0 \leq k \leq t} g_k$ encountered by the iterates during the algorithm after $t$ iterations satisfies: \vspace{-1mm}
\begin{equation} \label{eq:gBound}
\tilde{g}_t \leq \frac{\max\{2 h_0, C\} }{\sqrt{t+1}}  \quad \text{for $t \geq 0$} \, ,
\end{equation}
where $h_0 := f(\x^{(0)}) - \displaystyle \min_{\x \in \domain} f(\x)$ is the initial global suboptimality. It thus takes at most $O(\frac{1}{\epsilon^2})$ iterations to find an approximate stationary point with gap smaller than $\epsilon$.
\end{theorem}

The main idea of the proof is fairly simple and follows the spirit of the ones used for the gradient descent method. Basically, during FW, the objective $f$ is decreased by a quantity related to the gap $g_t$ at each iteration. As the maximum progress is bounded by the global minimum on~$\domain$ of~$f$, the gap~$g_t$ cannot always stay big, and the initial suboptimality $h_0 = f(\x^{(0)}) - \min_{\x \in \domain} f(\x)$ will control how big the gap can stay.

\begin{proof}
Let $\x_\stepsize
:=  \x^{(t)} + \stepsize \dd_t$ be the point obtained by moving with
step size $\stepsize$ in direction $\dd_t$, where $\dd_t := \st - \xt$ is the FW direction as defined 
by Algorithm~\ref{alg:FW}. By using $\s := \st$, $\x := \x^{(t)}$ and $\y := \x_\stepsize$ in the definition of the
curvature constant~$\Cf$~\eqref{eq:Cf}, and solving for $f(\x_\stepsize)$, we
get an affine invariant version of the standard descent lemma \citep[see e.g. (1.2.5) in][]{Nesterov:2004:lectures}:
\begin{equation} \label{eq:descentLemmaCf}
f(\x_\stepsize) \leq f(\x^{(t)}) + \stepsize \left\langle  \nabla
f(\x^{(t)}), \dd_t \right\rangle + \frac{\stepsize^2}{2} \Cf, \quad \text{valid $\forall
\stepsize \in [0,1]$}.
\end{equation}
Replacing the value of the FW gap $g_t$ in the above equation and substituting $C \geq \Cf$, we get:
\begin{equation} \label{eq:htStepsize}
f(\x_\stepsize) \leq f(\x^{(t)}) - \stepsize g_t + \frac{\stepsize^2}{2} C, \quad \text{valid $\forall
\stepsize \in [0,1]$}.
\end{equation}
We consider the best feasible step size that minimizes the quadratic upper bound on the RHS of~\eqref{eq:htStepsize}: $\stepsize^* := \min\{\frac{g_t}{C}, 1\}$. This is the same step size as used in option II of the algorithm ($f(\x^{(t+1)}) = f(\x_{\stepsize^*})$).  
In option I, the step size $\stepsize_t$ is obtained by line-search, and so
$f(\x^{(t+1)}) = f(\x_{\stepsize_t}) \leq  f(\x_\stepsize)$ $\forall
\stepsize \in [0,1]$ and thus $f(\x^{(t+1)}) \leq f(\x_{\stepsize^*})$. In both cases, we thus have:
\begin{equation} \label{eq:ht}
f(\x^{(t+1)}) \leq f(\x^{(t)}) - \min\left\{ \frac{{g_t}^2}{2 C} \, , \, g_t - \frac{C}{2} \indic_{ \{g_t > C \}} \right\} \, ,
\end{equation}
where $\indic_{ \{ \cdot \}}$ is an indicator function used to consider both possibilities of $\stepsize^* = \min\{\frac{g_t}{C}, 1\}$ in the same equation.
By recursively applying~\eqref{eq:ht}, we get:
\begin{equation} \label{eq:htPlusOne}
f(\x^{(t+1)}) \leq f(\x^{(0)}) - \sum_{k=0}^{t} \min\left\{ \frac{{g_k}^2}{2 C} \, , \, g_k - \frac{C}{2} \indic_{ \{g_k > C \}}\right\}.
\end{equation}
Now let $\tilde{g}_t := \min_{0\leq k \leq t} g_k$ be the minimal gap seen so far. Inequality~\eqref{eq:htPlusOne} then becomes:
\begin{equation} \label{eq:htInequality}
f(\x^{(t+1)})  \leq f(\x^{(0)}) - (t+1) \min\left\{ \frac{{\tilde{g}_t}^2}{2 C} \, ,\, \tilde{g}_t - \frac{C}{2} \indic_{ \{\tilde{g}_t > C \}} \right\}.
\end{equation}
We consider the two possibilities for the result of the $\min$. In both cases, we use the fact that $f(\x^{(0)}) - f(\x^{(t+1)}) \leq f(\x^{(0)}) - \min_{\x \in \domain} f(\x) =: h_0$ by definition and solve for $\tilde{g}_t$ in~\eqref{eq:htInequality}. In case that $\tilde{g}_t \leq C$, the first argument of the $\min$ is smaller and we get the claimed rate on $\tilde{g}_t$:
\begin{equation} \label{eq:gtRate}
\tilde{g}_t \leq \sqrt{ \frac{2 h_0 C}{t+1} } \,  .
\end{equation} 
In case that $\tilde{g}_t > C$ (in the first few iterations), we get that the initial condition $h_0$ is forgotten at a faster rate for $\tilde{g}_t$:
\begin{equation} \label{eq:gtRateInit}
\tilde{g}_t \leq  \frac{h_0}{t+1} + \frac{C}{2} \, .
\end{equation}
We note that this case is only relevant when $h_0 > C/2$; neither of the bounds~\eqref{eq:gtRate} and~\eqref{eq:gtRateInit} are then dominating each other as the inequality~\eqref{eq:gtRateInit} has a faster $O(1/t)$ rate but with the worse constant~$h_0$. We can also show that $\tilde{g}_t \leq C$ (for any $t \geq 0$) when $h_0 \leq C/2$. Indeed,
as we assumed that $\tilde{g}_t  > C$ to get~\eqref{eq:gtRateInit}, we have that~\eqref{eq:gtRateInit} then implies:
\begin{align}
C &< \frac{h_0}{t+1} + \frac{C}{2} \notag \\
t+1 &< \frac{2 h_0}{C} .  \label{eq:boundT}
\end{align}
If $h_0 \leq \frac{1}{2} C$, \eqref{eq:boundT} then yields a contradiction as $t \geq 0$, implying that $\tilde{g}_{t} \leq C$ for all $t\geq 0$ in this case.

From this analysis, we can summarize the bounds as:
\begin{equation} \label{eq:fullBound}
\tilde{g}_t \leq \left\{
\begin{aligned}
\frac{h_0}{t+1} + \frac{C}{2} \quad &\text{for} \quad t+1 \leq \frac{2 h_0}{C}   \, , \\
\sqrt{ \frac{2 h_0 C}{t+1}} \quad &\text{otherwise} \, .
\end{aligned} \right.
\end{equation}
We obtain the theorem statement by simplifying the first option in~\eqref{eq:fullBound} by using that it only happens when $h_0 > \frac{C}{2}$ and $t+1 \leq \frac{2 h_0}{C}$, and thus:
\begin{align*}
\frac{h_0}{t+1} + \frac{C}{2} &= \frac{h_0}{\sqrt{t+1}} \left( \frac{1}{\sqrt{t+1}} + \frac{C}{2h_0} \sqrt{t+1} \right) \\
&\leq \frac{h_0}{\sqrt{t+1}} \left( \frac{1}{\sqrt{t+1}} + \sqrt{\frac{C}{2h_0}} \right) \\
&\leq \frac{h_0}{\sqrt{t+1}} \left( \frac{1}{\sqrt{t+1}} + 1 \right) \\
&\leq \frac{2h_0} {\sqrt{t+1}} \quad \text{for $t\geq 0$}.
\end{align*}
By using $\sqrt{2h_0 C} \leq \max \{ 2 h_0, C \}$, we get the theorem statement.
\end{proof}

\section{Related work}

\paragraph{FW methods.} The only convergence rate for a FW-type algorithm on non-convex objectives that we are aware of is given in Theorem~7 of~\citet{yu2014GCG},\footnote{In this paper, they generalize the Frank-Wolfe algorithm to handle the (unconstrained) minimization of $f(\x) + g(\x)$, where $g$ is a non-smooth convex function for which $\min_{\x} \innerProd{\x}{\r} + g(\x)$ can be efficiently computed. The standard FW setup is recovered when $g$ is the characteristic function of a convex set~$\domain$, but they can also handle other types of $g$, such as the $\ell_1$ norm for example.} but they only cover non-adaptive step size versions (which does not apply to Algorithm~\ref{alg:FW}) and they can only obtain slower rates than~$O(1/\sqrt{t})$. 
\citet[Section 2.2]{bertsekas1999nonlinear} shows that any limit point of the sequence of iterates for the standard FW algorithm converges to a stationary point (though no rates are given). His proof only requires the function $f$ to be continuously differentiable.\footnote{We note that the gradient function being continuous on a compact set implies that it is uniformly continuous, with a \emph{modulus of continuity} function that characterizes its level of continuity. Different rates are obtained by assuming different levels of uniform continuity of the gradient function. The more standard one is assuming Lipschitz-continuity of the gradient, but other (slower) rates could be derived by using various growth levels of the modulus of continuity function.} He basically shows that the sequence of directions $\dd_t$ obtained by the FW algorithm is \emph{gradient related}, and then get the stationarity point convergence guarantees by using~\citep[Proposition 2.2.1]{bertsekas1999nonlinear}.
\citet[Note 5.5]{dunn1979rates} generalizes the standard rates for the FW method (in terms of global suboptimality) when running it on a class of quasi-convex functions of the form: $f(\x) := h(G(\x))$ where $h : \R \rightarrow\R$ is a strictly increasing real function with a continuous derivative, and $G : \R^d \rightarrow\R$ is a convex function with Lipschitz continuous gradient. These functions are quite special though: they are \emph{invex}, that is, all their stationarity points are also global optima.

\paragraph{Unconstrained gradient methods.} Our~$O(1/\sqrt{t})$ rate is analogous to the ones derived for projected gradient methods. In the unconstrained setting, \citet[Inequality~(1.2.15)]{Nesterov:2004:lectures} showed that the gradient descent method with line-search or a fixed step size of $\frac{1}{L}$, where $L$ is the Lipschitz constant of the gradient function, had the following convergence rate to a stationary point:
\begin{equation} \label{eq:gradientConvergence}
\min_{0 \leq k \leq t} \| \nabla f(\x^{(k)}) \| \leq \frac{\sqrt{2 h_0 L}}{\sqrt{t+1}} \quad \text{for $t \geq 0$} \, .
\end{equation}
We see that this rate is very similar to the one we give in~\eqref{eq:gtRate} in the proof of Theorem~\ref{thm:gConvergence}.
\citet{Cartis2010} also showed that the $O(1/\sqrt{t})$ rate was tight for the gradient descent method for an unconstrained objective. It is unclear though whether their example could be adapted to also show a lower bound for the FW method in the constrained setting, as their unidimensional example has a stationarity point only at~$+\infty$, which thus does not apply to a compact domain.

\paragraph{Constrained gradient methods.} In the constrained setting, several measures of non-stationarity have been considered for projected gradient methods. \citet{Cartis2012a} consider the \emph{first-order criticality measure} of~$f$ at~$\x$ which is similar to the FW gap~\eqref{eq:FWgap}, but replacing the maximization over $\domain$ in its definition to the more local $\domain \cap \mathbb{B}(\x)$, where $\mathbb{B}(\x)$ is the unit ball around~$\x$. This measure appears standard in the trust region method literature~\citep{Conn1993}. \citet{Cartis2012a} present an algorithm that gives a $O(1/\sqrt{t})$ rate on this measure. \citet{ghadimi2016nonconvexPG} considers instead the norm of the gradient mapping as a measure of non-stationarity.\footnote{The gradient mapping is defined for the more general proximal optimization setting, but we consider it here for the simple projected gradient setup. For a step size $\stepsize$, the gradient mapping is defined as $\frac{1}{\stepsize}(\x - \x^+)$ where $\x^+ := \text{Proj}_\domain(\x - \stepsize \nabla f(\x))$. If we let the $\stepsize \rightarrow 0$, then the gradient mapping becomes simply the negative of the projection of $-\nabla f(\x)$ on the solid tangent cone to~$\domain$ at~$\x$. When $\domain$ is the full space (unconstrained setting), then the gradient mapping becomes simply $\nabla f(\x)$. We use the ``gradient mapping'' terminology from~\citep[Definition 2.2.3]{Nesterov:2004:lectures} but with the notation from~\citep{ghadimi2016nonconvexPG}.} They show in~\citet[Corollary~1]{ghadimi2016nonconvexPG} that the simple projected gradient method with $\frac{1}{L}$ step size gives the same rate as given by~\eqref{eq:gradientConvergence} in the unconstrained setting, but using the norm of the gradient mapping on the LHS instead. They also later showed in~\citet{ghadimi2016accelerated} that the accelerated projected gradient method of Nesterov gave also the same $O(1/\sqrt{t})$ rate, but with a slightly better dependence on the Lipschitz constant.

\clearpage

\bibliographystyle{abbrvnat}
\bibliography{references}
\end{document}